\documentclass[a4paper,fleqn]{article} %Options in documentclass should be set to a4paper and fleqn.

\usepackage{modsim}
\usepackage{times}
\usepackage{natbib} %The three packages modsim, times and natbib are required.
\usepackage{amsmath, amssymb, amsthm} %Also recommend the standard AMS LaTeX maths packages.

\MODSIMhead{Y. Bokor {\it et al.} Stratified Space Learning} %This is the content of the headings in all pages (except the first page). The format should be author, title of paper. If the title is too long, then use ... at the end. If more than two authors, then please use the "et al." format with et al. in italic, for example A. Author {\it et al.}, Title of the paper.

\usepackage{rotating}
\usepackage{amsbsy,enumerate}
\usepackage{graphicx}
\usepackage{ccaption}
\usepackage{xcolor}
\usepackage[allbordercolors=white]{hyperref}
\newcommand{\orcid}[1] {\hspace*{-1.5mm} \href{https://orcid.org/#1}{\includegraphics[scale=0.3]{orcid1.jpg}}}

\newtheorem{definition}{Definition}

\makeatletter
\renewcommand{\fnum@figure}[1]{\textbf{\figurename~\thefigure}. }
\renewcommand{\fnum@table}[1]{\textbf{\tablename~\thetable}. }
\makeatother

\usepackage{pgfplots}
\usepackage{algorithmicx}
\usepackage{graphics}
\usepackage[noend]{algpseudocode}

\makeatletter

\def\BState{\State\hskip-\ALG@thistlm}

\makeatother

\begin{document}
% Defining Front matter.

\title{Stratified Space Learning: Reconstructing Embedded Graphs}

\author{\underline{Y. Bokor} \address[MSI]{\it{Mathematical Sciences Institue, Australian National University, Canberra, Australia}} %\orcid{0000-0002-4861-9174} 
and D. Grixti-Cheng \addressmark[MSI] and M. Hegland \addressmark[MSI] and S. Roberts \addressmark[MSI] and K. Turner \addressmark[MSI]} %underline the name of the presenting author. Use \addressmark[name-of-addressmark] after the name of an author who share the same affiliation (see modsim.tex for an example).

\email{yossi.bokor@anu.edu.au} %Email address of the presenting author only.

\date{July 2019}

%Keywords should be separated by commas and listed in Sentence case (first keyword with capital first letter and remaining keywords in lower case).

\begin{abstract}
	Many data-rich industries are interested in the efficient discovery and modelling of structures underlying large data sets, as it allows for the fast triage and dimension reduction of large volumes of data embedded in high dimensional spaces. The modelling of these underlying structures is also beneficial for the creation of simulated data that better represents real data. In particular, for systems testing in cases where the use of real data streams might prove impractical or otherwise undesirable. We seek to discover and model the structure by combining methods from topological data analysis with numerical modelling. As a first step in combining these two areas, we examine the recovery a linearly embedded graph $|G|$ given only a noisy point cloud sample  $X$ of $|G|$.

	An \emph{abstract graph} $G$ consists of two sets: a set of vertices $V$ and a set of edges $E$. An \emph{embedded graph} $|G|$ in $n$ dimensions is a geometric realisation of an abstract graph obtained by assigning to each vertex a unique coordinate vector in $\mathbb{R}^n$ and then each edge is the line segment between the coordinate vectors of the corresponding vertices. Given an embedded graph $|G| \subset \mathbb{R}^n$, a \emph{point cloud} sample of $|G|$ consists of a finite collection of points in $\mathbb{R}^n$ sampled from $|G|$, potentially with noisy perturbations. We will suppose that this sample has bounded noise of $\epsilon$ and is sufficiently sampled so every point in $|G|$ is within $\epsilon$ of some sample. Such a sample is called an $\epsilon$-sample.
	
	We can view this as a semiparametric model. Once the abstract graph is fixed we have a parametric model where the parameters are the locations of the vertices. In order to guarantee correctness of our algorithm, we will need to make some reasonable geometric assumptions on the embedded graph. 

	To learn the embedded graph $|G|$, we first learn the structure of the abstract graph $G$. We do this by assigning a dimension of either 0 or 1 to each $x \in X$ (depending on local topological structure) and then cluster points into groups representing embedded vertices or edges. Using local topological structure, we then assign to each abstract edge cluster a pair of abstract vertex clusters, to obtain the incidence relations of the abstract graph. Finally, we use nonlinear least squares regression to model the embedded graph $|G|$.

	The approach presented in this paper relies on topological concepts, such as stratification and local homology, which will be used in future research that will generalise this approach to embeddings on more general structures, in particular \emph{stratified spaces}.
	
	A stratified space $\mathbb{X}$ is a topological space with a partition into topological manifolds $\{X_i\}_{i \in I}$, called strata, such that for all $i$ and $j$, $X_i \cap X_j = \emptyset$, and if $X_i \cap \overline{X_j} \neq \emptyset$, then $X_i \subseteq \overline{X_j}$. 
	%We call this second condition the \emph{frontier axiom}.
	
	If we are interested in a specific topological structure, we can require that the strata satisfy aditional criteria, so that this structure is constant across each strata. In the case presented in this paper, the vertices of a graph (abstract or embedded) are the 0-strata, and the (open) edges are 1-strata.

	The authors are unaware of any algorithms which recover both the abstract graph, and model its embedding.

\end{abstract} %Abstract should NOT extend beyond the first page.

\maketitle

\section{INTRODUCTION}

    We seek to uncover information about an embedded $|G| \subset\mathbb{R}^n$ from a set $X\subset\mathbb{R}^n$ of noisy sample point $x^{(i)}$ which are close to $|G|$. We consider $|G|$ to be topological space, and with the topology induced by $\mathbb{R}^n$. Here we consider a \emph{stratified space} $|G|$ which consists of a partition $|G| = \bigcup_{i=1}^k |G|_i$ with strata $|G|_i$ overlapping at their boundaries.
    
    \begin{definition}
		A topologically stratified space $G$ is a space with a partition into topological manifolds $\{G_i\}_{i \in I}$, called strata, such that for all $i$ and $j$, $G_i \cap G_j = \emptyset$, and if $G_i \cap \overline{G_j} \neq \emptyset$, then $G_i \subseteq \overline{G_j}$. 
		%We call this second condition the \emph{frontier axiom}.
	\end{definition}
	
	One says that the strata $G_i$ are \emph{glued together} at their boundaries. The aim of this project is to recover the strata from the data points. Here the data points are an \emph{$\epsilon$-sample} $X \subset \mathbb{R}^n$ of $|G|$, which is a point cloud such that the \emph{Hausdorff} distance between $X$ and $|G|$ is bounded above by $\epsilon$. 
	Recall that the Hausdorff distance between two subsets $Y,Y' \subset \mathbb{R}^n$ is
	\begin{equation} d_H(Y,Y') := \max\{ \sup_{y \in Y} \inf_{y' \in Y'} d(y,y'), \sup_{y' \in Y'} \inf_{y \in Y} d(y,y') \},\end{equation}
	with $d(y,y')$ the standard Euclidean distance on $\mathbb{R}^n$.
	
	This is an inverse problem and the corresponding direct problem is generating data points from strata. The project consists of three logical steps:

	\begin{enumerate}
		\item the data is partitioned as \begin{equation} X = \bigcup_{j=1}^k X_j \end{equation} where the partition $X_j$ contain the points originating from stratum $G_j$,
		\item establish the dimension of the strata, and the glueing between them,
		\item concise numerical descriptions of the strata are computed.
  	\end{enumerate}

	When viewed as a semiparametric modelling problem, the first two steps determine the non-parametric components. For this we use topological data analysis. In contrast the last step numerically fits the parametric model discovered in the earlier procedures.
	
	In this paper we do not consider all stratified spaces and instead restrict to linear embeddings $|G|$ of abstract graphs $G$. An abstract graph $G$ is a stratified space with only 0- and 1dimensional strata, and each 1 dimensional strata (1-stratum) has two 0 dimension strata (0-stratum) as its boundary. We thus consider the 1-stratum $G_i^1$ to be a tuple of 0-strata $(G^0_{i_1}, G^0_{i_2})$. We often omit the superscript, in particular when writing 1-strata as a tuple. 

	%In addition to some topological and numerical background, we will discuss the special case where we have 0 dimensional strata (single point sets) and linear 1-dimensional strata (lines between two points). The abstract topology of the set $G$ is then defined by an \emph{(abstract) graph} with vertex set $V_0 = \{ i \mid \operatorname{dim}(G_i) = 0 \}$  and edge set $E_0$ where
	%  \begin{equation} (j_1,j_2) \in E_0 \Leftrightarrow
	%	  \text{there exists some $i\not\in V_0$ such that
	%  $G_i = \operatorname{conv}(g_{j_1},g_{j_2})$}. \end{equation}
	%Then we have the \emph{embedded graph} with vertex set $V = \{g_i \mid i \in V_0\}$ and edge set $E = \{\operatorname{conv}(g_{j_1},g_{j_2}) \mid (j_1,j_2) \in E_0\}$.
	
	To obtain a stratification of a graph, we will consider the \emph{local homology} of each strata. Given a topological space $G$, the local homology at a point $g$ describes the linear combinations of closed manifolds in a small region of $G$ around $g$. Examples of stratifications using other local topological structures can be found in \cite{stratifiedspaces}, \cite{Nanda2019}, \cite{sheafstratifications}, and \cite{localhomology}. \cite{stratifiedspaces} provide an overview of different types of stratifications, Nanda uses information about the local cohomology of cells in a regular CW complex to obtain a canonical stratification, \cite{sheafstratifications} apply sheaves to obtain stratifications of cell complexes, and \cite{localhomology} appeal to local homology to infer which points in a sample of a stratified space have been sampled from the same strata. Future research will expand the class of embedded stratified spaces we can recover from $\epsilon$-samples.
	
	\section{GEOMETRIC AND TOPOLOGICAL INTUITION}

	The first part of this semi-parametric problem is to identify the abstract graph structure and partition the points with respect to this abstract graph. To do this we classify samples topologically by their local neighbourhoods. This can be defined rigorously through the notion of local homology.
	Although the theoretical development of modifications of local homology to work with point cloud motivates and informs the algorithms here, for the sake of clarity and brevity we will focus on the computational geometry intuition behind the methods. This is possible as the local homology for graphs can be expresed in terms of connected components of apropriate subsets. For more complicated stratified spaces, more sophisticated machinery from algebraic topology will be required. A description relating this computational geometry to local homology is in the subsection below.
	
	We begin with the idea that whether a point on a graph is on an edge or a vertex can be determined by looking at how many times the graph intersects a sphere around that point of a small enough radius. For a point on an edge, the sphere will intersect twice, and for a degree $l$ vertex, it will intersect $l$ times.
	
	However, we do not have access to the underlying graph but instead only an $\epsilon$-sample. 
	Our sense of local has to be significantly larger than $\epsilon$ to be geometrically measurable. 
	This has multiple flow on effects in terms of how we analyse local environments to determine the abstract structure and how we partition the points into the vertices and edges.
	
	Since we have a sample of points that are not a subset of the underlying graph it is nonsensical to choose whether the samples lie \emph{on} a vertex or an edge. 
	Instead, we decide whether they are \emph{near} a vertex or an edge. We will classify points by the local geometry. For this we will fix a distance of $10\epsilon$ to be our sense of local. 
	That is for each sample point we look at all the samples within $10\epsilon$ of it and try to understand the local geometry at this $10\epsilon$ scale. 
	Does this geometric structure look that that of a neighbourhood of a vertex or of an edge? 
	This is the purpose of the dimension function within the algorithm.
	
	\begin{figure}[h!]
	
		\includegraphics[width=5cm]{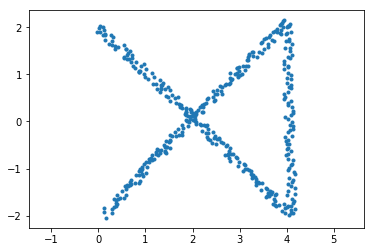}
		\includegraphics[width=5cm]{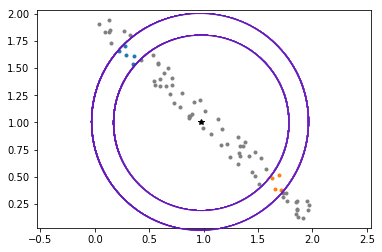}
		\includegraphics[width=5cm]{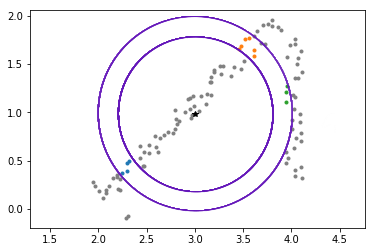}\\
		\, (a) \hspace{4.5cm} (b) \hspace{4.5cm} (c)\\
		\includegraphics[width=5cm]{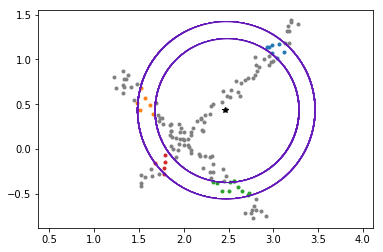}
		\includegraphics[width=5cm]{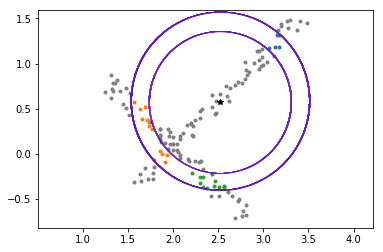}
		\includegraphics[width=5cm]{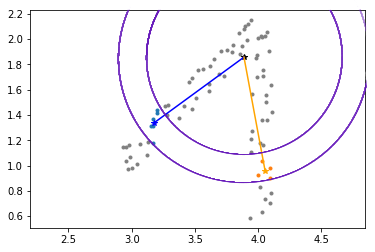}
		\, (d) \hspace{4.5cm} (e) \hspace{4.5cm} (f)\\
		\caption{(a) contains the $\epsilon$-sample. (b)-(f) illustrate diiferent cases of local neighbourhood. The colours indicate different connected components in the annular region.}\label{geo}
	\end{figure}
	
	When we had the underlying graph we could measure the unique point on each nearby edge which intersected the sphere. 
	With an $\epsilon$-sample we have to use an annulus instead of a sphere in order to ensure that there exists a sample point for each intersection. 
	By making the difference of the inner and outer radius to be $2\epsilon$ (as in the region with distance between $8\epsilon$ and $10\epsilon$ from the sample) we can guarantee that a ball of radius $\epsilon$ around a location in the graph is contained in the annular region and hence there will be a sample point in the $\epsilon$-sample. 
	There can be multiple samples in this annular region for each edge and so we consider equivalence classes of points that would be samples of the that part of the edge. 
	This can be achieved by considering connected components where we connect two points when they are within $3\epsilon$.    
	
	There are two cases where we declare that a sample is near an edge (or rather that it is not near a vertex). 
	The typical case is that the samples lying in an $10\epsilon$ ball around sample $q$ are noisy samples of a straight line that passes at most $\epsilon$ from $q$. 
	This set of points are all connected at threshold $2\epsilon$ and there are two connected components at threshold $3\epsilon$ for the set of samples restricted to the annulus. 
	Furthermore these connected components lie close to directly opposite each other. 
	The angle between the centroids is at least $2\arccos(1/4)$. 
	An example is shown in Figure \ref{geo} (b). 
	The other case occurs when we are closer to a vertex where there is an acute angle between adjacent edges. 
	Here the samples lying in an $10\epsilon$ ball are not all connected at threshold $2\epsilon$. 
	An example is shown in Figure \ref{geo} (c). 
	
	Given a lower bound of $\pi/6$ for the angle between adjacent vertices, we can guarantee one of these two cases will occur for every sample at least $15\epsilon$ away from every vertex. Thus, each such sample will be classified as near an edge by the dimension function.
	
	We also want to show that the dimension function declares all samples sufficiently close to a vertex are of dimension 0. 
	Near a vertex of degree larger than $2$, we will want the number of connected components to be at least $3$ for all samples sufficiently close to the vertex. This will depend on a lower bound on the angle between adjacent edges and for our algorithm (with the inner radius set at $8\epsilon$) a sufficient condition is a lower bound of $\pi/6$. With such a bound, any sample $q$ with $\|q-v\|<2\epsilon$ the samples within an $10\epsilon$ ball of $q$ are all connected at a threshold of $2\epsilon$ and each of the edges will contribute a different connected component within the annulus centred at $q$. 
	An example is illustrated in Figure \ref{geo} (f).
	
	Stronger statements can be made in terms of the existence of connected components representing each of a suitable subsets of edges near a sample point. 
	This covers other cases such as in Figure \ref{geo} (d). 
	We can guarantee that the samples declared to be of dimension 0 spread far enough to separate the points declared of dimension 1 on the different edges, so that they appear as separate clusters in $\text{dim}^{-1}(1)$.   
	
	To detect vertices of degree $2$ we have to introduce an angle condition. When the annulus has two connected components we measure the angle between the vectors pointing to the center of these two components. If we lie near an edge then the angle should be at least $2\arccos(1/4)$. In contrast if the angle at the vertex is $\theta$ and $q$ is within $4\epsilon$ of both edges then the computed angle between the connected components is bounded from above by $\theta+\arcsin(1/4) +\arcsin(5/8)$. Since our algorithm declares the dimension is $0$ once the angle is below $2\arccos(1/4)$ we can guarantee identifying the abstract graph structure when $\theta\leq \pi/2$. This covers other cases such as in Figure \ref{geo} (e).
	
	The dimension function effectively partitions the samples into those assigned to a vertex and those assigned to an edge we can then further partition the samples into the connected components with appropriate thresholds. 
	We use $3\epsilon$ as the threshold between samples in $\text{dim}^{-1}(1)$ to say they are of the same edge.
	A low threshold is needed to ensure that we do not merge adjacent edges. 
	We use the much larger $10\epsilon$ when looking to see with samples in $\text{dim}^{-1}(0)$ are connected. 
	We can use this higher threshold because of the geometric assumptions about minimum distances between vertices.
	
	Once we have a partition of the samples into the clusters for $\text{dim}^{-1}(0)$ and $\text{dim}^{-1}(1)$ we can read off the abstract graph. The $\text{dim}^{-1}(0)$ and $\text{dim}^{-1}(1)$ clusters are in bijection with the vertices and the edges respectively. Furthermore, for each edge $e\in E$ with corresponding $\text{dim}^{-1}(1)$ cluster $c(e)$ there are exactly two $\text{dim}^{-1}(0)$ clusters with points within $3\epsilon$ of $c(e)$. These correspond to the vertices that bound $e$ in the abstract graph.
	
	We then wish to model the locations of the vertices as the vertices completely determine the embedding of the graph. We do this using a least squares estimator which will minimise the sums of distances from the samples to the embedded graph where this distance is the distance to a vertex if the sample is in $\text{dim}^{-1}(0)$ and the distance is the distance to the corresponding edge if the sample is in $\text{dim}^{-1}(1)$. It should be noted that this is naturally a biased estimator which is visable in the simulations. The bias is for vertices to drift in the direction where there a more edges. A future direction of research is to reduce, or  estimate and correct for this bias.
	
	\subsection{Relationship to local homology}
	For $Y$ a topological space, the homology groups $H_i(Y)$ are algrabic summaries of the topological structure of $Y$. We have a homology group for each dimension $i$: $H_0(Y)$ describing the connected components of $Y$ and $H_1(Y)$ the space of loops. E.g. The homology of the wedge of $N$ circles is $\text{rank}(H_1(\bigvee_N S^1))=N$. 
	
	To measure the local topological structure at a point $q$ we can construct a quotient space where everything far from $q$ is quotiented to a single point. The local homology at $q$ is defined at the limit of the homology of $Y/\{y\in Y:d(y,q)>r\}$ as $r$ tends to $0$. The the case where the $Y$ is an embedded graph then the local homology of a vertex with degree $k$ is that of the wedge of $k-1$ copies of $S^1$. The local homology of a point along an edge is that of $\bigvee_{k-1} S^1$. 
	
	Since we do not have access to the underlying graph but instead an $\epsilon$-sample it is not meaningful to take the limit as $r$ goes to $0$. Instead the practice in topological data analysis is to use the sample point cloud to estimate the $R$-local homology which uses a fixed radius $r=R$ in the construction of the relevant quotient spaces. In the case of an embedded graph $Y$ with all vertices separated by at least $2R$ and $q$ a point near $X$, this $R$-local homology is that of the wedge of $k-1$ spheres where $k$ is the number of points in $Y\cap \partial B(q,R)$. Note that if an edge intersects twice if and only if it is its own separate connected components in $X\cap B(q,R)$. By counting the connected components in $Y\cap B(q,R)$ and $Y\cap \partial B(q,R)$ we can completely determine the topological structure of embedded graph $Y$ near $q$. 
	
\section{ALGORITHM}

	The steps 2 through 5 of the algorithm use topological data analysis while the last step is a nonlinear least-squares fit.

	\subsubsection{Step 1: Construct $\mathfrak{G}_{3\epsilon}(X)$}
		This step constructs an embedded graph $\mathfrak{G}_{3\epsilon}(X)$ on the sample $X$, by connecting points $p,q \in X$ if $d(p,q) \le 3 \epsilon$.

	\subsubsection{Step 2: Dimension Function}
		Given a point $q \in X$, we construct the subgraph of $\mathfrak{G}_{3\epsilon}(X)$ generated by the points within $10 \epsilon$ of $q$, and restrict the edges . If this graph is disconnected, the function outputs dimension 1. Otherwise, we consider the subcomplex generated by the points in the annulus of radii $8 \epsilon$ and $10 \epsilon$ around $q$. If the number of connected components is not 2, we assign dimension $0$, if it is 2, we apply the angle test.

		\begin{enumerate}
			\item Take the subgraph $\mathfrak{G}_q$ of $\mathfrak{G}_{3\epsilon}(X)$ generated by points $p \in X$ such that $d(p,q) \le 10 \epsilon$.
			\item Remove edges between points $p, p'$ if $d(p,p') > 2 \epsilon$.
			\item For $p,p'$ vertices in $\mathfrak{G}_q$, add an edge between $p$ and $p'$ if $d(p,p') \le 2 \epsilon$.
			\item If the number of connected components in $\mathfrak{G}_q$ is not $1$, output dimension $1$.
			\item Else, consider subgraph $\mathfrak{G}_q^{(1)}$ of $\mathfrak{G}_{3\epsilon}(X)$ generated by $p \in \mathfrak{G}_q$ such that $d(p,q) \ge 8 \epsilon$.
			\item If the number of connected components in $\mathfrak{G}_q^{(1)}$ is not 2, return dimension 0.
			\item Else, do Step 2a.
		\end{enumerate}
		%{\color{red} im not sure how we want to write the pseudo code, so I will just write it roughly in enumerate}

	\subsubsection{Step 2a: Angle test}
		For points $p$ with 2 connected components in the annulus with radii ($8\epsilon, 10\epsilon)$, we calculate the angle between the components by taking the mid-point of each cluster, and calculating the angle between the lines joining $p$ and each mid-point. If this angle is less that $2 \arccos(1/4)$, the dimension function outputs 0, otherwise it outputs 1.

		\begin{enumerate}
			\item Find average of coordinates of points in the two connected components.
			\item Calculate angle between the line segments from averages to $q$.
			\item If angle is less that $2 \arccos(1/4)$ return dimension 0.
			\item Else return dimension 1.
		\end{enumerate}

	\subsubsection{Step 3: Finding the vertices}
		Having determined the dimension of each point in the sample, we consider the subcomplex generated by the points we have assigned dimension 0. Considering the connected components of this subcomplex, we obtain the clusters of points surroudning vertices.

		\begin{enumerate}
			\item Generate subgraph of $\mathfrak{G}_{3 \epsilon}(X)$ on the vertices with dimension function 0.
			\item Map each connected component of the subgraph to a vertex of an abstract graph.
		\end{enumerate}

	\subsubsection{Step 4: Finding the edges}
		Having determined the dimension of each point in the sample, we consider the subcomplex generated by the points we have assigned dimension 1. Considering the connected components of this subcomplex, we obtain the clusters of points surrounding the edges.

		\begin{enumerate}
			\item Generate subgraph of $\mathfrak{G}_{3 \epsilon}(X)$ on the edges with dimension function 1.
			\item Map each connected component of the subgraph to a vertex of an abstract graph.
		\end{enumerate}

	\subsubsection{Step 5: Assigning boundary vertices to edges}
		For each edge cluster, we find the two vertex clusters which have points connected to those in the edge cluster in the full complex.

		\begin{enumerate}
			\item For each connected component $G^1_i$ in the preimage of 1 under the dimension function, we add in points in the preimage of 0 which are within $3\epsilon$cite of $G^1_i$. 
			\item Add an edge to the abstract graph $(citei,j)$ if points from the $i^{th}$ and $j^{th}$ connected components in preimage of $0$ are added above. cite
		\end{enumerate}

\subsubsection{Step 6: Nonlinear least squares fitting}

This is the final step in recovering the embedded graph $|G|$. Thus far, we have the abstract graph $(V_0, E_0)$ of $G$. 
We also determined for each data point $x^{(i)}$ the stratum $G_j$ which it is closest to. 
We will now determine the geometric information used to embed
$X$ in $\mathbb{R}^{d}$. For the 0-dimensional strata $X_j$ this is just one point $x_i \in \mathbb{R}^n$.
The set of 0 dimensional strata is just $\{x_j \mid j \in V_0\}$.
The 1-dimensional strata are of the form $X_j = \{\theta x_{j_1} + (1-\theta) x_{j_2}\mid \theta \in (0,1)\}$
and the set of 1-dimensional strata is $\{ \theta x_{j_1} + (1-\theta) x_{j_2} \mid (j_1,j_2)\in E,
\theta \in (0,1)\}.$ Thus we need to determine the $x_j$ for $j\in V_0$. We will use a least-squares approach
where $x_j$ are such that they minimise the sum of Euclidean distances squared between data points and their
associated strata. For the 0 dimensional strata we introduce a partial objective as
	\begin{equation} \phi_i(x_1,\ldots,x_{k_v},\theta_i) = \lVert p^{(i)} - x_{j(i)}\rVert^2.\end{equation}
Here $\theta_i = 0$.
For 1-dimensional strata we include a parameter $\theta_i$ which models the local coordinate of the
data point in $X_j$ and get a partial objective
\begin{equation}
\phi_i(x_1,\ldots,x_{k_v}, \theta_i) =
  \lVert p^{(i)} - \theta_i x_{j_1(i)} - (1-\theta_i) x_{j_2(i)}\rVert^2.
	\end{equation}
We then get an optimisation problem with objective
\begin{equation} \Phi(x_1,\ldots,x_{k_v}, \theta_1,\ldots,\theta_n) = \sum_{i=1}^n \phi_i(x_1,\ldots,x_{k_v}, \theta_i).\end{equation}
We minimise this $\Phi$ and have the following conditions for $\theta_i$:
\begin{itemize}
	\item $\theta_i \in [0,1]$
  \item $\theta_i = 0$ if $j(i) \in V_0$.
\end{itemize}

This optimisation problem is then solved using a code from scipy.optimize. It implements the
trust-region reflexive method and uses a Gauss-Newton approach. For more details and
references see \cite{scipy}. As the problem is sparse and converges within a finite
number of steps in practice, the complexity of this approach is approximately $O(n)$.

\section{SIMULATIONS}

\begin{figure}[h!]

	\includegraphics[width=4.8cm]{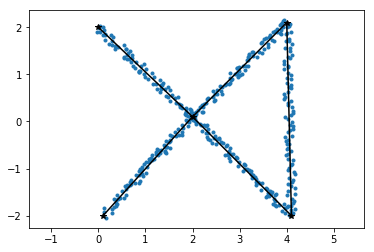}
	\includegraphics[width=4.8cm]{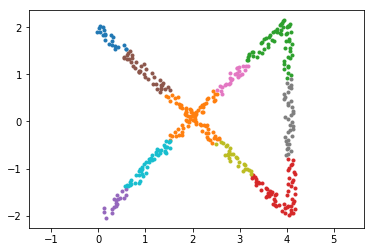}
	\includegraphics[width=4.8cm]{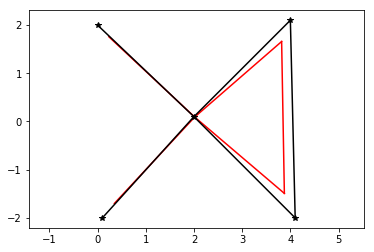}
	\caption{Left: An embedding $|G|_2$ of $G$ in $\mathbb{R}^2$ (black) and an $\epsilon$-sample $X_2$ of $|G|_2$ (red). Middle: partioniong of $X_2$ obtained from algorithm. Right: The embedded graph $|G|_2$ (black) and the model $\tilde{|G|_2}$ (red).} \label{2d}
\end{figure}

\begin{figure}[h!]
 	\includegraphics[width=4.8cm]{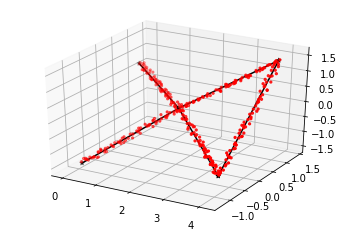}
 	\includegraphics[width=4.8cm]{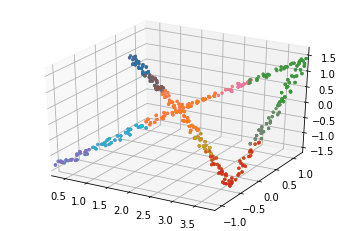} 
 	\includegraphics[width=4.8cm]{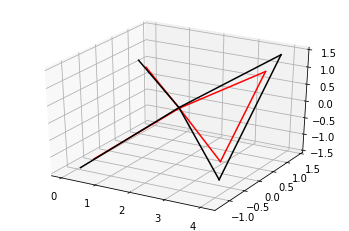}
 	\caption{Left: An embedding $|G|_3$ of $G$ in $\mathbb{R}^3$ (black) and an $\epsilon$-sample $X_3$ of $|G|_3$ (red). Middle: partioniong of $X_3$ obtained from algorithm. Right: The embedded graph $|G|_3$ (black) and the model $\tilde{|G|_3}$ (red).} \label{3d}
 \end{figure}

	In this section, let $G$ be the abstract graph with vertex set $\{0,1,2,3,4\}$ and edge set $\{(0,1), (1,2), (2,3), (3,1)\}$ and $\epsilon =  0.1$. We consider two embeddings of $G$: $ |G|_2 \subset \mathbb{R}^2$ and  $ |G|_3 \subset \mathbb{R}^3$. We will follow the same procedure for both. Taking an $\epsilon$-sample $X$, we will partition $X$, and then use these partitions to model the embedding. We used \cite{Hunter2007} for generating the images. See Figure \ref{2d} and Figure \ref{3d}. We provide these figures for illustrative purposes, and the graphs were chosen to display the cases discussed earlier. We are unaware of any algorithms which both model the embedding $|G|$ and recover the abstract structure of $G$ from an $\varepsilon$-sample $X$. 

\section{CONCLUSIONS}
 We have provided an algorithm that can learn embedded graphs under reasonable geometric constraints. This semiparametric modelling algorithm requires as input an $\epsilon$-sample of the embeded graph, and consists of two main components. First, it learns the abstract graph structure (the non-parametric component), and then using non-linear least squares regression, it finds an optimal embedding of vertices (the parameters of the model). Future directions include reduction, estimation and correction of the bias, expanding the algorithm to piecewise linear embedded structures, and modifications to deal with semi-algebraic sets.

\section*{Acknowledgement}
	The first author was supported by an Australian Government Research Training Program (RTP) Stipend and RTP Fee-Offset Scholarship through the Australian National University. The first author would also like to thank James Morgan for various helpful discussions.  %{\color{red} Any others we want/need to acknowledge?}
	
\bibliography{stratified_space_learning}
\bibliographystyle{chicago} %use chicago style of referencing.
\end{document}